\magnification=1200
\input amssym.def

 \input amssym.tex
 \font\newrm =cmr10 at 24pt
\def\bul{\raise .9pt\hbox{\newrm .\kern-.105em } }
\def\fr{\frak}

 \baselineskip 20pt
 
 \def\h{\hbox{ }}
 
 \def\p{{\fr p}}
 
 \def\r{{\fr r}}
 
 \def\m{{\fr m}}
 \def\n{{\fr n}}
 \def\a{{\fr a}}
 \def\d{{\fr d}}
 
 \def\ss{{\fr s}}
 
 \def\b{{\fr b}}
 \def\cc{{\fr c}}
 \def\hh{{\fr h}}
 
 \def\ee{{\fr e}}

 \def\g{{\fr g}}

 \def\q{{\fr q}}

 \def\<{\le}
 \def\>{\ge}

 \def\s{{\h\subset\h}}

 \def\mapright#1
  {\smash{\mathop
  {\longrightarrow}
  \limits^{#1}}}
 \def\mapdown#1{\Big\downarrow
  \rlap{$\vcenter
  {\hbox{$#1$}}$}}
  
 \def\kk#1{{\kern .4 em} #1}

\font\twelverm=cmr12 at 14pt
\hsize = 31pc
\vsize = 45pc
\overfullrule = 0pt
	\font\authorfont=cmr9
\font\ninerm=cmr9
\centerline{\twelverm Dirac Cohomology for the Cubic Dirac Operator} 
\vskip 1.5pc
\baselineskip=11pt
\vskip8pt
\centerline{\authorfont BERTRAM KOSTANT\footnote*{\ninerm
Research supported in part by NSF grant
DMS-9625941 and in part by the \hfil\break KG\&G Foundation}}\vskip 2pc
\baselineskip 20pt
\noindent {\bf Abstract.} Let $\g$ be a complex semisimple Lie algebra and let $\r\s \g$
be any  reductive Lie subalgebra such that $B|\r$ is
nonsingular where $B$ is the Killing form of $\g$. Let $Z(r)$ and
$Z(\g)$ be, respectively, the centers of the enveloping algebras of 
$\r$ and $\g$. Using a Harish-Chandra isomorphism one has a
homomorphism $\eta:Z(\g)\to Z(\r)$ which, by a well-known result
of H. Cartan, yields the the relative Lie algebra cohomology
$H(\g,\r)$. 

Let
$V$ be any $\g$- module. For the case where $\r$ is a
symmetric subalgebra, Vogan has defined the Dirac cohomology
$Dir(V)$ of
$V$. Using the cubic Dirac operator we extend his definition to the
case where $\r$ is arbitrary subject to the condition stated above. 
We then generalize results of Huang-Pand$\check {\hbox{z}}$i\'c
on a proof of a conjecture of Vogan. In particular $Dir(V)$ has a
 structure of a $Z(\r)$-module relative to a ``diagonal"
homomorphism $\gamma:Z(\r)\to End\,Dir(V)$. In case $V$ admits an
infinitesimal character $\chi$ and $I$ is the identity operator on
Dir(V) we prove $$\gamma\circ \eta = \chi\,\,I\eqno (A)$$ In
addition we also prove that $V$ always exists (in fact $V$ can taken
to be an object in Category O) such that $Dir(V)\neq 0$. If $\r$ has
the same rank as $\g$ and
$V$ is irreducible and finite dimensional, then (A) generalizes a
result of Gross-Kostant-Ramond-Sternberg.

\vskip 1.5pc

\centerline{\bf 0. Introduction}\vskip 1.5pc 0.1. Let
$\g$ be a complex semisimple Lie algebra and let $(x,y)$ be a nonsingular symmetric
invariant bilinear form $B_{\g}$ on $\g$. Let $\r\s
\g$ be any reductive Lie subalgebra of $\g$ such that $B_{\r} = B_{\g}\vert \r$ is
nonsingular. Let
$\p$ be the $B_{\g}$-orthocomplement of $\r$ in $\g$ so that $[\r,\p]\s \p$ and one has the
direct sum $\g = \r+ \p$. Let $B_{\p}= B_{\g}\vert \p$ so that $B_{\p}$ is nonsingular
and let $C(\p)$ be the Clifford algebra over $\p$ with respect to $B_{\p}$. Then there exists
a homomorphism $\nu_*:\r\to C(\p)$ such that $[x,y] = [\nu_*(x),y]$ for $x\in \r$ and
$y\in \p$ where the bracket on the right side is commutation in $C(\p)$. See \S 1.5 in
[K1]. One then has a homomorphism $$\zeta:U(\r)\to U(\g)\otimes C(\p)\eqno (0.1)$$ so that
$\zeta(x) = x\otimes 1 + 1\otimes \nu_*(x)$ for $x\in \r$. This defines the structure of an
$\r$-module on $U(\g)\otimes C(\p)$. We have defined an element 
$\square\in U(\g)\otimes C(\p)$ in [K1] and have referred to $\square$ as the cubic Dirac
operator. The definition of $\square$ is recalled in \S 2.2 of the present paper. In case
$\r$ is a symmetric subalgebra of
$\g$ the cubic term in
$\square$ vanishes and $\square$ is the more familiar Dirac operator under consideration in
[HP] and [P]. The result main result in [HP] works in the general case considered here and
one has a unique homomorphism
$$\eta_{\r}:Z(\g)\to Z(\r)$$ so that for any $p\in Z(\g)$ one has $$p\otimes 1-
\zeta(\eta_r(p)) = \square\,\omega + \omega \square \eqno (0.2)$$ for some $\omega\in
(U(\g)\otimes C^{odd}(\p))^{\r}$. See Appendix. We will determine the homomorphism
$\eta_{\r}$ in the generality under consideration here. 

A subspace $\ss\s \g$ will be said to be $(\r,\p)$ split if $\ss = \ss\cap \r +\ss\cap\p$.
In such a case we will write $\ss_{\r}= \ss\cap \r $ and $\ss_{\p} = \ss\cap\p$. Let $\b$
be a Borel subalgebra of $\g$. Let $\hh\s \b$ be a Cartan subalgebra of $\g$ and let $\n$
be the nilradical of $\b$. Then $\b$ and $\hh$ can (and will) be chosen so that they are
both $(\r,\p)$ split. Then $\b_{\r}$ is a Borel subalgebra of $\r$ and $\hh_{\r}$ is a
Cartan subalgebra of $\r$. It follows also that $\n$ is $(\r,\p)$-split. The subspace
$\n_{\p}\s \p$ is isotropic with respect to $B_{\p}$. Let $u\in C(\p)$ be the product, in
any order, of a basis of $\n_{\p}$. Let $L= C(\p)\,u$ so that $L$ is a left ideal in
$C(\p)$ and is in particular a $C(\p)$-module with respect to left multiplication. Now
assume that $V$ is a $U(\g)$-module. The separate actions of $U(\g)$ and $C(\p)$ define a
homomorphism $$\xi_{V}:U(\g)\otimes C(\p)\to End\,(V\otimes L)$$ Extending Vogan's definition
to the present case one defines the Dirac cohomology $\break$ $H_D(V\otimes L)$ so that 
$$H_D(V\otimes L) = Ker\,\square_{V}/(Ker\,\square_{V}\cap Im\,\square_{V})\eqno (0.3)$$
where
$\square_V= \xi_{V}(\square)$. The action of $\xi_V(Z(\g)\otimes 1)$ on $V\otimes L$
defines a $Z(\g)$-module structure on $H_D(V\otimes L)$. Also if $\zeta_V=
\xi_V\circ \zeta$ then the action of $\zeta_V(Z(\r))$ on $V\otimes L$ defines a
$Z(\r)$-module structure on 
$H_D(V\otimes L)$. As a consequence of (0.2) one has, for any $p\in Z(\g)$,
$$\xi_V(p\otimes 1) = \zeta_V(\eta_{\r}(p))\,\,\hbox{on}\,\,H_D(V\otimes L)\eqno (0.4)$$ This
raises the question as to whether or not $H_D(V\otimes L)$ vanishes. 

Let $\lambda\in \hh^*$ and let $V_{\lambda}$ be the unique irreducible (Category O)
$U(\g)$-module with highest weight $\lambda$ (with respect to $\b$). Let $0\neq
v_{\lambda}\in V_{\lambda}$ be a corresponding highest weight vector. If $dim\,\hh_{\p}= k$
then $\Bbb C\,v_{\lambda}\otimes C(\hh_{\p})u$ is a $2^k$-dimensional subspace of
$V_{\lambda}\otimes L$. Here $C(\hh_{\p})\s C(\p)$ is the Clifford algebra over $\hh_{\p}$.
Let $\rho\in \hh^*$ have its usual meaning. In this paper we will prove \vskip 1pc {\bf
Theorem 0.1.} {\it Choose $\lambda$ so that $\lambda + \rho$ vanishes on $\hh_{\p}$ (e.g.,
$\lambda = -\rho$ if $\r=0$). Then 
$\Bbb C\,v_{\lambda}\otimes C(\hh_{\p})u \s Ker\,\square_{V_{\lambda}}$. Furthermore the
map from cocycle to cohomology defines an injection $$\Bbb C\,v_{\lambda}\otimes
C(\hh_{\p})u \to H_D(V_{\lambda}\otimes L)\eqno (0.5)$$ In particular
$H_D(V_{\lambda}\otimes L)\neq 0$.} \vskip 1pc Since there is no restriction on
$\lambda\vert\hh_{\r}$ we can compute $\eta_{\r}$. 

Let $\phi_o:\hh\to \hh_{\r}$ be the
projection relative to the decomposition $\hh = \hh_{\r} + \hh_{\p}$. Then $\phi_o$
extends to a homomorphism $S(\hh)\to S(\hh_{\r})$ and clearly induces a homomorphism
$$\phi:S(\hh)^{W_{\g}}\to S(\hh_{\r})^{W_{\r}}$$ where $W_{\g}$ and $W_{\r}$ are the
respective Weyl groups of $\hh$ relative to $\g$ and $\hh_{\r}$ relative to $\r$. Given the
fact that $\r$ is essentially an arbitrary reductive Lie subalgebra of $\g$ the following
result established here is a strong generalization of Theorem 5.5 in [HP]. As will be noted
in \S 5 in this paper it is also a generalization of Proposition 3.43, (5.18) and (5.19) in
[K1]. 
\vskip 1pc {\bf Theorem 0.2.} {\it The map
$\eta_{\r}:Z(\g)\to Z(\r)$ is uniquely determined so that the following diagram is
commutative. In the diagram the vertical maps are the Harish-Chandra isomorphisms.
$$
\def\mapright#1{\smash{
\mathop{\longrightarrow}\limits^{#1}}}
\def\mapdown#1{\Big\downarrow
\rlap{$\vcenter{\hbox{$\scriptstyle#1$}}$}}
\matrix
{Z(\g)&\mapright {\eta_{\r}}\, &Z(\r)\cr
\mapdown {A_{\g}}&&\mapdown {A_{\r}}\cr
S(\hh)^{W_{\g}}&\mapright{\phi}&S(\hh_{\r})^{W_{\r}}\cr}$$} \vskip 1pc
The map $\phi$ is well known in the theory of the cohomology of compact homogeneous
spaces. Actually what is utilized in that theory is the map $S(\hh^*)^{W_{\g}}\to
S(\hh_{\r}^*)^{W_{\r}}$ induced by restriction of functions. However this is same as
$\phi$ if $\hh$ and $\hh^*$ are identified and $\hh_{\r}$ and $\hh_{\r}^*$ are
identified using $B_{\g}$. Assume $G$ is a compact connected semisimple Lie group and $\g$ is
the complexification of $Lie\,G$. Let $R\s G$ be any connected compact subgroup and let $\r$
be the complexification of $Lie\,R$. Obviously we can choose $B_{\g}$ so that
$B_{\g}\vert \r$ is nonsingular (e.g., let $B_{\g}$ be the Killing form). The 
map $\eta_{\r}$ induces the structure of a $Z(\g)$-module on $Z(\r)$. On the other
hand the infinitesimal character for the module $V_{\lambda}$ when $\lambda =-\rho$
defines the structure of a $Z(\g)$-module on $\Bbb C$. As a consequence of a
well-known theorem of H. Cartan (see \S 9 in [C]) one has \vskip 1pc {\bf Theorem 0.3.} {\it
There exists an isomorphism $$H^*(G/R,\Bbb C) \cong Tor_*^{Z(\g)}(\Bbb C, Z(\r))\eqno
(0.6)$$}\vskip 1pc In \S 5 we reformulate certain results in [K1] using Dirac cohomology.
\vskip 1pc 0.2. We wish to thank David Vogan for many profitable conversations and for
introducing us to his Dirac cohomology concept in the case where $\r$ is a symmetric
subalgebra of $\g$. We also wish to acknowledge the strong impact made upon us by the main
result in [HP].
\vskip 1.5pc

	\centerline{\bf 1. Preliminaries}\vskip 1.5pc
1.1. Let $\g$ be a semisimple complex Lie algebra and let $B_{\g}$ be a
nonsingular $ad$-invariant symmetric bilinear form $(x,y)$ on
$\g$. Let $\r\s\g$ be a reductive Lie subalgebra and assume that $B_{\r}=
B_{\g}\vert \r$ is nonsingular. Let $\p$ be the $B_{\g}$-orthocomplement of 
$\r$ in $\g$ and let $B_{\p} = B_{\g}\vert \p$. Then of course $$\g = \r + \p$$
and $[\r,\p]\s \p$. Let $\hh_{\r}$ be a Cartan subalgebra of $\r$ and let $\hh\s\g$
be a Cartan subalgebra of $\g$ containing $\hh_{\r}$. Of course $B_{\r}\vert
\hh_r$ and $B_{\g}\vert \hh$ are nonsingular. Let $\hh_{\p}$ be the
orthocomplement of $\hh_{\r}$ in $\hh$ so that $$\hh = \hh_{\r} +\hh_{\p}\eqno (1.1)$$ and
$$\eqalign{\hh_{\r} &= \hh\cap \r\cr 
\hh_{\p} &=\hh\cap \p\cr}$$ Let $\Delta\s \hh^*$ be the set of roots for
$(\hh,\g)$ and for each $\varphi\in \Delta$ let $e_{\varphi}\in \g$ be a
corresponding root vector. We normalize the choice so that
$(e_{\varphi},e_{-\varphi})= 1$. Let
$\g^0$ be the centralizer of $\hh_{\r}$ in $\g$ and let $\Delta^0 = \{\varphi\in
\Delta\mid \varphi(x) = 0,\,\forall x\in \hh_{\r}\}$ so that
$$\g^0 = \hh +\sum_{\varphi\in \Delta^0} \Bbb C\,e_{\varphi}\eqno (1.2)$$ Let
$\hh^{\#}$ be the real space of hyperbolic elements in $\hh$ and let $\kappa:\hh\to
\hh^{\#}$ be the real projection which vanishes on $i\hh^{\#}$. Since
$\hh_{\r}$ is complex there clearly exists
$f_{\r}\in \kappa(\hh_{\r})$ such that if $\g^{f_{\r}}$ is the centralizer of $f_{\r}$
in $\g$ then $\g^{f_{\r}} = \g^0$. But $f_{\r}$ defines a parabolic Lie subalgebra
of $\q$ of $\g$ where $\g^0$ is a Levi factor of $\q$ and the $nilrad\,\q$ is the
span of all eigenvectors of $ad\,f_{\r}$ with positive eigenalues. Clearly
$ad\,f_{\r}$ stabilizes both $\r$ and $\p$ and hence $$nilrad\,\,\q = \n_{\r} +
\p^+\eqno (1.3)$$ where $\n_{\r} = \r\cap nilrad\,\,\q$ and $\p^+ = \p\cap nilrad\,\,\q$.
Since clearly $$\g^0\cap \r = \hh_{\r}\eqno (1.4)$$ it follows that $$\b_{\r} =
\hh_{\r} + \n_{\r}\eqno (1.5)$$ is a Borel subalgebra of $\r$ and $\n_{\r}$ is the
nilradical of
$\b_{\r}$. Furthermore (1.4) implies that $$\g^0 = \hh_{\r}  + \p^0\eqno (1.6)$$ where
$\p^0= \g^0\cap \p$, is
an orthogonal decomposition with respect to the (obviously) nonsingular bilinear
form $B_{\g}\vert \g^0$. Let $\cc =Cent\,\g^0$. Since of course $\hh_{\r}\s \cc$ one
has
$$\cc = \hh_{\r} + \cc_{\p}\eqno (1.7)$$ where $\cc_{\p}=\cc\cap \p$, is an orthogonal
decomposition with respect to the (obviously) nonsingular bilinear form
$B_{\g}\vert \cc$. Of course $\cc\s\hh$ so that $\cc_{\p}\s \hh^{\p}$. Let
$\d_{\p}$ be the orthocomplement of $\cc_{\p}$ in $\hh_{\p}$ so that $$\hh_{\p} =
\cc_{\p} +
\d_{\p}\eqno (1.8)$$ is an orthogonal decomposition. 
\vskip 1pc {\bf Remark 1.1} Note that (1.2) and (1.4) imply that $$\p^0 = \hh_{\p}
+\sum_{\varphi\in
\Delta^0}\Bbb C e_{\varphi}\eqno (1.9)$$ and that $\p^0$ is a reductive Lie subalgebra
of $\g$ which happens to lie in $\p$. Furthermore $\hh_{\p}$ is a Cartan subalgebra
of $\p^0$ and
$\cc_{\p}\s \hh_{\p}$ is the center of $\p^0$. In particular
$$\p^0 = \cc_{\p} + [\p^0,\p^0]\eqno (1.10)$$ is an orthogonal decomposition
and (1.8) implies that $\d_{\p}$ is a Cartan subalgebra of the semisimple Lie algebra
$[\p^0,\p^0]$. Obviously $$[\p^0,\p^0] = \d_{\p} +
\sum_{\varphi\in
\Delta^0}\Bbb C e_{\varphi}\eqno (1.11)$$ and (1.11) is the decomposition of
$[\p^0,\p^0]$ as the sum of a Cartan subalgebra and corresponding root
spaces. Let $\p'$ be the orthocomplement of $\p^0$ in $\p$. Clearly $\p'$ is stable
under $ad\,\hh_{\r}$ and hence $\p'$ is stable under $ad\,f_{\r}$. Let $\p^-$ be
the span of all eigenvectors of $ad\,f_{\r}$ in $\p'$ with negative eigenvalues.
Obviously $\p' = \p^+ + \p^-$ so that one has direct sums $$\eqalign{\p&= \p^0
+\p'\cr &=\p^0 + \p^+ +\p^-\cr}\eqno (1.12)$$ Let $\Gamma\s \hh_{\r}^*$ be the set of
all weights for the adjoint action of $\hh_{\r}$ on $\p'$ and for any $\mu\in
\Gamma$ let $\p^{\mu}\s\p'$ be the corresponding weight space so that one has the
direct sum $$\p' = \sum_{\mu\in \Gamma}\p^{\mu}\eqno (1.13)$$ It is clear that any
$\mu\in
\Gamma$ extends uniquely to a linear functional (to be identified with $\mu$) on the complex
subspace of $\hh$ spanned by $\kappa(\hh_{\r})$. One has a partition
$\Gamma = \Gamma_+ \cup \Gamma_-$ where
$$\Gamma_+\,(\hbox{resp.}\,\Gamma_-)\,=\{\mu\in \Gamma\mid \mu(f_{\r})
>0\,(\hbox{resp.}\,\mu(f_{\r})<0)\}$$\vskip .5pc {\bf Remark 1.2.} By considering the
action of $ad\,\hh_r$ a standard argument implies that for $\mu,\nu\in \Gamma$
one has that $\p^{\mu}$ is $B_{\p}$-orthogonal to $\p^{\nu}$ if $\nu\neq -\mu$. But
since $B_{\p}\vert\p'$ is clearly nonsingular one has that $\Gamma = -\Gamma$ and
$\p^{\mu}$ is nonsingularly paired to $\p^{-\mu}$ for any $\mu\in \Gamma$. It then
follows that $\Gamma_- = -\Gamma_+$ and $$\eqalign{\p^+ &= \sum_{\mu\in
\Gamma_+}\p^{\mu}\cr \p^- &= \sum_{\mu\in
\Gamma_+}\p^{-\mu}\cr}\eqno (1.14)$$ Obviously there exists a closed Weyl chamber
$C\s \hh^{\#}$ such that $f_{\r}\in C$. Let $f\in\hh^{\#}$ be an element in
the interior of $C$ so that, in particular, $f\in \hh$ is regular and hyperbolic.
One defines a choice of positive roots $\Delta_+\s \Delta$ by putting
$\Delta_+=\{\varphi\in \Delta\mid \varphi(f) >0\}$. Let $\Delta_- = -\Delta_-$. Let $\b\s\g$
be the Borel subalgebra defined by putting $$\b = \hh + \sum_{\varphi\in \Delta_+} \Bbb
C\,e_{\varphi}$$ Let $\n = [\b,\b]$ be the nilradical of $\b$. Since $f_{\r}\in C
$ one readily has $$\eqalign{\b&\s \q\cr nilrad\,\,\q&\s \n\cr}\eqno (1.15)$$ Let $\n^0
= \n\cap \p^0$, $\Delta^0_+ = \Delta_+\cap \Delta^0$ and $\Delta^0_- = -\Delta^0_+$ so that
$$\n^0 = 
\sum_{\varphi\in \Delta_+^0} \Bbb C\,e_{\varphi}\eqno (1.16)$$ It then follows
from (1.3),(1.5) and (1.15) that $$\n = \n_{\r} + \n^0 + \p^+\eqno (1.17)$$ noting now that
$\n_{\r}$, the nilradical of the Borel subalgebra $\b_{\r}$ of $\r$, is given by
$$\n_r = \n\cap\r\eqno (1.18)$$ and also $$\n^0 + \p^+ = \n\cap \p\eqno (1.19)$$\vskip
1.5pc \centerline{\bf 2. Dirac cocycles}\vskip 1.5pc 2.1. Let $C(\p)$ be the
Clifford algebra over $\p$ with respect to $B_{\p}$. As in \S 1.5 of [K1] we
identify the underlying linear spaces of $C(\p)$ and the exterior algebra
$\wedge\,\p$ and understand that there are two multiplications in $C(\p)$. If
$w,z\in C(\p)$ then $w\,z$ denotes the Clifford product and $w\wedge z$ the
exterior product of $w$ and $z$. If $w\in
\wedge^k\p$ and $z\in \wedge^{k'}\p$ then one knows $$w\,z-w\wedge z\in
\sum_{j=0}^{k+k'-2}\wedge^j\p\eqno (2.1)$$ (for an argument see e.g., (1.6) in [K1]).
The bilinear form
$B_{\p}$ on
$\p$ extends to an nonsingular bilinear form $(w,z)$ on $C(\p)$, to be denoted by $B_{C(\p)}$,
so that if
$w\in
\wedge^k\p$ and $z\in \wedge^{k'}\p$ then $(w,z) = 0$ if $k\neq k'$. If $k=k'$ then 
$(w,z) = det\,(w_i,z_j)$ where $w=w_1\wedge\cdots\wedge w_k$ and $z
=z_1\wedge\cdots\wedge z_k$ for $z_i,w_j\in \p$. It
is immediate then that $$\m_{\p} = \n^0 + \p^+\eqno (2.2)$$ is a $B_{\p}$-isotropic subspace of
$\p$. However since $\m_{\p}$ is $B_{\p}$-isotropic it follows that Clifford product and
exterior product are the same for elements in $\m_{\p}$. Let $u_0$ be the product of all the
root vectors $e_{\varphi}$ for $\varphi\in
\Delta^0_+$ in some order and let $u_{+}$ be the product of a basis of $\p^+$ in some
order. Put $u = u_0\,u_+$ so that, in $C(\p)$, $$ z\,u = 0\,\,\,\,\,\forall z\in
\m_{p}\eqno (2.3)$$ Let $L\s C(\p)$ be the left ideal $$L = C(\p)\,u\eqno (2.4)$$
In particular $L$ is a $C(\p)$-module under left multiplication. Let $C(\hh_{\r})$ be
the Clifford algebra over $\hh_{\r}$ so that $C(\hh_{\r})$ is a subalgebra of
$C(\p)$. \vskip 1pc {\bf Proposition 2.1.} {\it The map $$C(\hh_{\r})\to L, \qquad
a\mapsto au$$ is injective. Furthermore for any $z \in \m_{\p}$ and $a\in
C(\hh_{\r})$ one has $$z\,a\,u = 0\eqno (2.5)$$}\vskip 1pc {\bf Proof.} The first
statement is a consequence of (2.1) and the fact that $\hh_{\r}\cap \m_{\p} = 0$.
The equation (2.5) follows from (2.3) and the fact that $\hh_{\r}$ is
$B_{\p}$-orthogonal to $\m_{\p}$. QED \vskip .5pc 2.2. Let $U(\a)$ be the universal
enveloping algebra of $\a$ where $\a\s\g$ is any Lie subalgebra. We are mainly concerned 
here with the algebra tensor product $U(\g)\otimes C(\p)$. If $x,y\in\p$ then no confusion
should arise from $x\otimes y$ as an element in $U(\g)\otimes C(\p)$. The left
factor $x$ is taken to be in $U(\g)$ and the right factor $y$ is taken to be in
$C(\p)$. In \S 2.1 of [K1] we introduced an element $\square\in U(\g)\otimes C(\p)$ which we
referred to as a cubic Dirac operator (see \S 0.23 in [K1]). We recall the definition of
$\square$. Let $I$ be an index set having cardinality equal to $dim\,\p$. Then $\square =
\square' + \square''$ where if $\{z_i\},\,i\in I$, is an orthonormal basis of $\p$ with
respect to $B_{\p}$ one has $$\square' = \sum_{i\in I} z_i\otimes z_i\eqno (2.6)$$ and
$$\square'' = 1\otimes v\eqno (2.7)$$ where $v\in \wedge^3\p$ is such that for any
$x,x',x''\in \p$ one has $$([x,x'],x'') = -2(v,x\wedge x '\wedge x'')\eqno (2.8)$$ See (1.20)
in [K1]. 

Now let $\lambda\in \hh^*$ be arbitrary and let $V_{\lambda}$ be the unique irreducible
highest module for $U(\g)$ with highest weight (relative to $\b$) $\lambda$. We recall that
$\lambda$ extends uniquely to a character on $\b$, $z\mapsto \lambda(z)$, which
necessarily vanishes on $\n$, and if $\Bbb C_{\lambda}$ is the corresponding 1-dimensional
$U(\b)$-module then $V_{\lambda}$ is the quotient of the Verma module
$U(\g)\otimes_{U(\b)}\Bbb C_{\lambda}$ by the unique maximal proper submodule. Let
$0\neq v_{\lambda}\in V_{\lambda}$ be a highest weight vector so that $z\,v_{\lambda} =
\lambda(z)v_{\lambda}$ for any $z\in \b$. Now let $V_{\lambda,L} = V_{\lambda}\otimes L$ so
the action of $U(\g)$ on $V_{\lambda}$ and $C(\p)$ on $L$ defines an algebra homomorphism
$$\xi_{\lambda}: U(\g)\otimes C(\p) \to End\,V_{\lambda,L}\eqno (2.9)$$ Let $a\in
C(\hh_{\r}$ and put $$v_{\lambda,a} = v_{\lambda}\otimes a\,u\eqno (2.10)$$ The element
$v_{\lambda,a}\in V_{\lambda,L}$ is nonzero, by Proposition 2.1, if $a\neq 0$. Our principal
goal now is to compute $\xi_{\lambda}(\square)v_{\lambda,a}$. 

For any $\nu\in \hh^*$ let $z_{\nu}\in \hh$ be the element corresponding to $\nu$ with
respect to the isomorphism $\hh^*\to \hh$ defined by $B_{\g}\vert \hh$. Thus $(z,z_{\nu}) =
\nu(z)$ for any $z\in \hh$. But now, by (1.1), there uniquely exists $x_{\nu}\in \hh_{\r}$ and
$y_{\nu}\in \hh_{\p}$ such that $$z_{\nu} = x_{\nu} + y_{\nu}\eqno (2.11)$$

We will first deal with $\xi_{\lambda}(\square')v_{\lambda,a}$. Let $\{b_i\},\,i\in I$, be
any basis of $\p$ and let $\{d_i\},\,i\in I, $ be the dual basis with respect to $B_{\p}$.
It is clear from (2.6) that $\square'$ can be rewritten as the sum $$\square' =\sum_{i\in
I}b_i\otimes d_i\eqno (2.12)$$ We will now partition the index set $I$, first, as a union of
three parts $$ I = I^{\hh}\cup I^0\cup I'\eqno (2.13)$$ where $\{b_j\},\,j\in I^{\hh}$, is an
orthonormal basis of $\hh_{\p}$. Next $\{b_k\},\,k\in I^0,  = \{e_{\varphi}\},\,\varphi\in
\Delta^0$, and $\{b_{m}\},\,m\in I'$, is a basis of $\p'$. We further refine the choice of
the basis by partitioning $$I' = \bigcup_{\mu\in \Gamma}I^{\mu}\eqno (2.14)$$ so that
$\{b_m\},\,m\in I^{\mu}$, is a basis of $\p^{\mu}$. By Remark 1.2 we can make the choice so
that
$$\{b_{m'}\},\,m'\in I^{-\mu}, = \{d_m\},\,m\in I^{\mu},\eqno (2.15)$$ Note also that if $j\in
I^0$ and $b_j= e_{\varphi}$ for $\varphi\in \Delta^0$ then necessarily one has $$d_j=
e_{-\varphi}\eqno (2.16)$$ In addition for any $j\in I^{\hh}$ one then clearly has $$d_j =
b_j\eqno (2.17)$$ \vskip .5pc {\bf Remark 2.2.} Note that, by (2.15), (2.16) and (2.17),
setwise
$$\{b_i\},\,j\in I, = \{d_i\},\,i\in I,$$ Also for any $i\in I^0\cup I'$ one has, by (2.2),
$$\hbox{either}\,\,b_i\in \m_{\p}\,\,\hbox{or}\,\,d_i\in \m_{\p}\eqno (2.18)$$\vskip .5pc
{\bf Lemma 2.3.} {\it For any $\lambda\in \hh^*$ and $a\in C(\hh_{\p})$ one has
$$\xi_{\lambda}(\square')v_{\lambda,a} = v_{\lambda}\otimes y_{\lambda}\,a\,u\eqno
(2.19)$$}\vskip 1pc {\bf Proof.} By (1.17) and (2.2) one has $$\m_{\p}\s \n\eqno (2.20)$$ It
then follows from (2.5) and (2.18) that, for any $i\in I^0\cup I'$, $$\xi_{\lambda}(b_i\otimes
d_i)v_{\lambda,a} =0\eqno (2.21)$$ However for $j\in I^{\hh}$ one clearly has, by (2.5),
$$\xi_{\lambda}(b_j\otimes d_j)v_{\lambda,a} =\lambda(b_j)v_{\lambda}\otimes b_j\,a\,u\eqno
(2.22)$$ But clearly $\sum_{j\in I^{\hh}}\lambda(b_j)\,b_j = y_{\lambda}$. This proves
(2.19). QED\vskip 1pc 2.3. We will now compute $\xi_{\lambda}(\square'')v_{\lambda,a}$. To do
so we first introduce a simple ordering in $I$. We will choose the ordering so that $i<j $ if
$i\in I^{\hh}$ and
$j\in I^0$, and also if $i\in I^0$ and
$j\in I'$. Also we fix the order so that if $i,j\in I^0$ and $b_i = e_{\varphi},\,b_j = 
e_{\varphi'}$ for some $\varphi,\varphi'\in \Delta^0$ then $i< j $ if $\varphi\in \Delta^0_+$
and
$\varphi'\in
\Delta^0_-$. In addition if $i,j\in I'$ then $i<j$ if $b_i\in \p^+$ and $b_j\in \p_-$. Let
$T'$ be the set of all ordered triples $\{i,j,k\}$ where $i,j,k\in I$ and $i<j<k$. The set
$\{b_i\wedge b_j\wedge b_k\},\,\{i,j,k\}\in T'$, is of course a basis of $\wedge^3\p$. The
dual basis with respect to $B_{C(\p)}\vert \wedge^3\p$ is clearly $\{d_i\wedge d_j\wedge
d_k\},\,
\{i,j,k\}\in T'$. But now we can write $$v= \sum_{\{i,j,k\}\in T'} c_{ijk}\,\,b_i\wedge
b_j\wedge b_k\eqno (2.23)$$ for $c_{ijk}\in \Bbb C$. But then for any $\{i,j,k\}\in T'$ one has
$$([d_i,d_j],d_k) = -2\, c_{ijk}\eqno (2.24)$$ by (2.8). 

But it is clear from our choice of basis that, for any $i\in I$, $b_i$ is a weight
vector for some weight $\gamma_i\in \hh_{r}^*$ with respect to the action of $ad\,\hh_{r}$ on
$\p$. Also it is clear that
$\gamma_i\in \Gamma\cup \{0\}$. Note then it follows from Remark 2.2 that $d_i$ is a weight
vector with weight $-\gamma_i$. But then it follows from (2.24) that $c_{ijk} \neq 0$ implies
that
$\gamma_i + \gamma_j + \gamma_k = 0$. Thus if $T = \{\{i,j,k\}\in T'\mid \gamma_i + \gamma_j +
\gamma_k = 0\}$ then one has $$v = \sum_{\{i,j,k\}\in T} c_{ijk}\,\,b_i\wedge b_j\wedge
b_k\eqno (2.25)$$ Let $i\in I$. Then obviously $$\eqalign{\gamma_i = 0\,\,&\iff b_i \in \p^0 
\cr &\iff i\in I^{\hh}\cup I^0\cr
\gamma_i\neq 0\,\,&\iff b_i\in \p'\cr
&\iff i\in I'\cr}\eqno (2.26)$$ Now let $T_0 = \{\{i,j,k\}\in T\mid \gamma_i = \gamma_j =
\gamma_k = 0\}$ so that if $\{i,j,k\}\in T$ then $\{i,j,k\}\in T_0$ if and only if
$\{b_i,b_j,b_k\}\s \p^0$. Let $T_1$ be the complement of $T_0$ in $T$. For $\varepsilon =0,1$,
let 
$$v^{(\varepsilon)} =
\sum_{\{i,j,k\}\in T_{\varepsilon}} c_{ijk}\,\,b_i\wedge b_j\wedge b_k$$ so that $$v =
v^{(0)} + v^{(1)}\eqno (2.27)$$ and $$\square'' = \square_0'' + \square_1''\eqno (2.28)$$
where $\square_{\varepsilon}'' = 1\otimes v^{(\varepsilon)}$. Let $\rho_0 = {1\over 2}
\sum_{\varphi\in \Delta^0_+} \varphi$. \vskip 1pc {\bf Lemma 2.4.} {\it For any $\lambda\in
\hh^*$ and $a\in C(\hh_{\p})$ one has
$$\xi_{\lambda}(\square_0'')v_{\lambda,a} = v_{\lambda}\otimes y_{\rho_0}\,a\,u\eqno
(2.29)$$}\vskip 1pc {\bf Proof.}. Let $\{i,j,k\}\in T_0$. Now if $i\in I^0$ then $b_i =
e_{\varphi_1}$ for $\varphi_1\in \Delta^0$. But then by the order relation on $I$ one has
$b_j= e_{\varphi_2}$ and $b_k= e_{\varphi_3}$ where also $\{\varphi_2,\varphi_3\}\s
\Delta_0$. But then by (2.16) and (2.24) one has $c_{ijk} = 0$ if $\varphi_1 + \varphi_2
+\varphi_3\neq 0$. However if $\varphi_1 + \varphi_2 +\varphi_3 =  0$ then by the ordering one
has $\varphi_1\in \Delta^0_+$. But then $b_i\wedge b_j\wedge b_k = b_j\,b_k\,b_i$ since the 
elements $b_i,b_j$ and $b_k$ are mutually orthogonal so that exterior and Clifford
multiplication are the same. But then $\xi_{\lambda}(b_i\wedge b_j\wedge b_k)v_{\lambda,a}=0$
by (2.5). Thus in computing the left side of (2.29) we can ignore all the terms in the
definition of $v^{(0)}$ for which $i\in I^0$. Now assume that $\{i,j,k\}\in T_0$ and $i\in
I^{\hh}$. But then if also $j\in I^{\hh}$ one has $c_{ijk}= 0 $ by (2.17) and (2.24). If
$j\in I^0$ then also $k\in I^0$ so that for some $\varphi,\varphi'\in \Delta^0$ one has $b_j =
e_{\varphi}$ and $b_k= e_{\varphi'}$. But then $c_{ijk} \neq 0$ implies $\varphi +\varphi'=
0$ by (2.16) and (2.24). But then $\varphi\in \Delta^0_+$ and $\varphi' = -\varphi$.
Moreover (2.16) and (2.17) imply $$c_{ijk} = {1\over 2}\,\varphi(b_i)\eqno (2.30)$$ But since
$b_i$ is orthogonal to $e_{\beta}$ for any $\beta\in \Delta$ one has
$$v^{(0)}\,a\,u= \sum_{i\in I^{\hh},\,\varphi\in \Delta^0_+} {1\over 2}\,\varphi(b_i)\,\, b_i
(e_{\varphi}\wedge e_{-\varphi})\,a\,u\,\eqno (2.31)$$ But now if $z,w\in \p $ and $(z,w)=1$
then $$z\wedge w = -w\,z +1 \eqno (2.32)$$ by e.g., (1.6) in [K1]. Thus $e_{\varphi}\wedge
e_{-\varphi}= e_{\varphi}\,e_{\varphi} +1$ in (2.31). However $e_{\varphi}\,a\,u = 0$ by
(2.5). Consequently $$\eqalign{v^{(0)}\,a\,u &= (\sum_{i\in I^{\hh},\,\varphi\in \Delta^0_+}
{1\over 2}\,\varphi(b_i)\,\, b_i)\,a\,u\cr 
&=(\sum_{i\in I^{\hh}}
\,\rho_0(b_i)\,\, b_i)\,a\,u\cr
&= y_{\rho_0}\,a\,u\cr}\eqno (2.33)$$ But of course (2.33) implies (2.29). QED\vskip 1pc
2.4. We will now determine $v^{(1)}\,a\,u$. We have assumed that $a\in C(\hh_{\p})$. However
for later purposes we want $a$ to be more general. Let $C(\p^0)\s C(\p) $ be the Clifford
algebra over
$\p^0$ with respect to $B_{\p}\vert \p^0$. The argument establishing (2.5) also establishes
$$z\,a\,u= 0\eqno (2.34)$$ for $z\in \p^+$ and $a\in C(\p^0)$. 
\vskip 1pc {\bf Lemma 2.5.} {\it Assume $a\in C(\p^0)$. Let
$\{i,j,k\}\in T_1$. Then if $(b_i\wedge b_j\wedge b_k)\,a\,u\neq 0$ one has $b_i\in \p^0$
(i.e., $i\in I^\hh\cup I^0$), $j\in I^{\mu}$, for some $\mu\in \Gamma_+$ and $b_k = d_j$.
Furthermore in such a case 
$$(b_i\wedge b_j\wedge b_k)\,a\,u = b_i\,a\,u$$ Moreover in this case $$c_{ijk}
= {1\over 2}([d_i, b_j],d_j)\eqno (2.35)$$}
\vskip 1pc {\bf Proof.} One must have $j\in I'$, since otherwise $b_j\in \p^0$ in which case
$b_i\in
\p^0$ (by the ordering). However this implies $b_k\in \p^0$ since $$\gamma(i) + \gamma(j)
+\gamma(k) =0\eqno (2.36)$$ But then
$\{i,j,k\}\in T_0$ which is a contradiction. Thus $j\in I'$ and hence $k\in I'$. But then if
$i\in I'$ the relation (2.35) and the ordering implies that $b_i\in \p^+$ and $b_i,b_j$ and
$b_k$ are mutually orthognal so that $b_i\wedge b_j\wedge b_k = b_j\,b_k\,b_i$. But
$b_i\,a\,u= 0 $ by (2.34). Hence our nonvanishing assumption implies that $b_i\in
\p^0$. But now $j\in I^{\mu}$, by (2.14), for some $\mu\in \Gamma$. But then $k\in I^{-\mu}$
by (2.36). Hence $\mu\in \Gamma_+$ by (1.14) and the ordering in $I$. In particular
$b_j\in \p^+$. Thus the nonvanishing assumption in the lemma implies $$b_i\wedge b_j\wedge b_k
= b_i( b_j\wedge b_k)$$ But if $b_j$ is orthogonal to $b_k$ then $b_j\wedge b_k = -b_k\,b_j$.
But
$b_j\,a\,u = 0$ by (2.34). In the remaining case where $b_k$ is not orthogonal to
$b_j$ one has $b_k = d_j$ by (2.15). However $b_j\wedge d_j = -d_jb_j +1$ by (2.32). But
$b_j\,a\,u= 0$ by (2.34). With the exception of (2.35) this proves the lemma. But
the equality $b_k = d_j$ implies $b_j = d_k$ again by (2.15). But then (2.24) implies
that $c_{ijk} = -{1\over 2}([d_i, d_j],b_j)$. But, by the invariance of $B_{\g}$, $([d_i,
d_j],b_j) = -([d_i, b_j],d_j)$. This proves (2.35). QED
\vskip 1pc Now let
$\Delta^1_+ =
\{\varphi\in
\Delta_+\mid
\varphi\vert\hh_{\r}\in
\Gamma_+\}$ and let $\Delta^2_+ = \{\varphi\in \Delta_+\mid \varphi\vert \hh_{\r}\notin
\Gamma_+\cup {0}\}$. Then one has a partition $$\Delta_+ = \Delta_+^0\cup \Delta_+^1\cup
\Delta^2_0\eqno (2.37)$$ We have already defined $\rho_0$. Let $\rho, \rho_1$ and $\rho_2$ be
defined similarly where $\Delta_+,\Delta_+^1$ and $\Delta_+^2$ respectively replace
$\Delta_+^0$. Thus $$\rho = \rho_0 + \rho_1 + \rho_2\eqno (2.38)$$ Hence $$y_{\rho} =
y_{\rho_0}+y_{\rho_1} +y_{\rho_2}\eqno (2.39)$$ \vskip .5pc {\bf Lemma 2.6.} {\it One has
$$y_{\rho_2}= 0\eqno (2.40)$$} \vskip 1pc {\bf Proof.} Let $\varphi\in \Delta^2_+$ and let
$e_{\varphi}^{\r}\in \r$ and $e_{\varphi}^{\p}\in \p$ be such that $e_{\varphi}
=e_{\varphi}^{\r} + e_{\varphi}^{\p}$. But since $ad\,\hh_{\r}$ stabilizes both $\r$ and $\p$
it follows that $e_{\varphi}^{\r}$ and $e_{\varphi}^{\p}$ are $ad\,\hh_{\r}$-weight vectors
for the weight $\varphi\vert \hh_{\r}$. But, by definition, $\varphi\vert \hh_{\r}\notin
\Gamma_+\cup
\{0\}$. But since, clearly, $\varphi(h_{\r})\geq 0$ one has $e_{\varphi}^{\p}= 0$.
Hence $e_{\varphi}\in \r$. But $[y,e_{\varphi}] = \varphi(y)\,e_{\varphi}\in \r$ for any $y\in
\hh_{\p}$. But $[\r,\p]\s \p$. Thus $$\varphi\vert\hh_{\p} = 0,\qquad \forall \varphi\in
\Delta_+^2\eqno (2.41)$$ This proves (2.40). QED\vskip 1pc Now for any $\mu\in \Gamma$ let
$\Delta^1_{\mu} = \{\varphi\in \Delta\mid \varphi\vert \hh_{\r}=\mu\}$. Also for any $\mu\in
\Gamma$ let
$\g^{\mu}$ be the weight space in $\g$ corresponding to the weight $\mu$ with respect to the 
action of $ad\,\hh_r$ on $\g$. Clearly $\g^{\mu}$ and $\g^{-\mu}$ are nonsingularly paired
by $B_{\g}$. In fact $\g^{\mu}$ is clearly stable under the adjoint action of
$\g^0$ and, since $\hh\s \g^0$, it follows immediately that $\{e_{\varphi}\},\,\varphi \in
\Delta^1_{\mu}$, is a basis of $\g^{\mu}$ and $\{e_{-\varphi}\},\,\varphi \in
\Delta^1_{\mu}$, is the $B_{\g}$-dual basis in $\g^{-\mu}$. It follows therefore that, for any
 $z\in \g^0$, $${1\over 2}tr\,ad\,z\vert \g^{\mu} = {1\over 2}\sum_{\varphi\in
\Delta^1_{\mu}}([z,e_{\varphi}],e_{-\varphi})\eqno (2.42)$$ In particular if $\rho_1^{\mu} =
{1\over 2}\, \sum_{\varphi\in \Delta^1_{\mu}}\varphi$ and if $z\in \hh_{\p}$, then $${1\over
2}\,tr\,ad\,z\vert \g^{\mu}= \rho_1^{\mu}(z)\eqno (2.43)$$ But clearly $\sum_{\mu\in
\Gamma_+}\rho_1^{\mu} = \rho_1$ since one readily has the partition $$\Delta_+^1 =
\bigcup_{\mu\in \Gamma_+}\Delta^1_{\mu}$$ Thus for $z\in \hh_{\p}$, $${1\over
2}\sum_{\mu\in \Gamma_+}\,tr\,ad\,z\vert \g^{\mu}= \rho_1(z)\eqno (2.44)$$ For any
$\mu\in \Gamma$ let $\r^{\mu} =\r\cap \g^{\mu}$. \vskip 1pc {\bf Lemma 2.7.} {\it For any
$\mu\in \Gamma$ one has $$\g^{\mu} = \r^{\mu} + \p^{\mu}\eqno (2.45)$$} \vskip 1pc {\bf
Proof.} Obviously the right side of (2.45) is contained in the left side. Conversely let $w\in
\g^{\mu}$ and let
$w_{\r}\in
\r$ and
$w_{\p}$ be such that $w= w_{\r} +w_{\p}$. But clearly both components $w_{\r}$ and $w_{\p}$
are also weight vectors of $ad\,\hh_{r}$ with weight $\mu$. Thus $w_{\p}\in \p^{\mu}$ and
obviously $w_{\r}\in \r^{\mu}$. QED \vskip 1pc Let $\mu\in \Gamma_+$. Obviously $\r^{\mu}$
is $B_{\g}$-nonsingularly paired to $r^{-\mu}$. We already know the same is true if
$\p^{\mu}$ and $\p^{-\mu}$ replace $\r^{\mu}$ and $\r^{-\mu}$. If $z\in \p^0\s \g^0$, to
compute $tr\,ad\,z\vert \g^{\mu}$ instead of using a basis of root vectors as we did in
(2.42) we can use the basis $\{b_j\},\,j\in I^{\mu},$ of $\p^{\mu}$ together with some
basis of $\r^{\mu}$. But since $[\p,\r]\s \p$ it follows that $[z,\r^{\mu}]\s\p^{\mu}$.
Thus we need consider only the basis $\{b_j\},\,j\in I^{\mu},$ of $\p^{\mu}$ to compute
the trace. Thus for any $i\in I^\hh\cup I^0$ (so that $b_i,d_i\in \p^0$) and $\mu\in \Gamma_+$
one has $$tr\,ad\,\,d_i\vert \g^{\mu} = \sum_{j\in I^{\mu}}([d_i, b_j],d_j)\eqno (2.46)$$ But
then, by Lemma 2.5, for $a\in C(\p^0)$ one has $${1\over 2}\,\sum_{\mu\in
\Gamma_+}\,(\sum_{i\in I^{\hh}\cup I^0}(tr\,ad\,\,d_i\vert
\g^{\mu})b_i)\,a\,u = v^{(1)}\,a\,u\eqno (2.47)$$ But $d_i\in [\g^0,\g^0]$ for $i\in I^0$ by
(1.11) and (2.16). Hence $tr\,ad\,\,d_i\vert \g^{\mu} = 0$ for $i\in I^0$ and any $\mu\in
\Gamma_+$. Thus, by (2.15), (2.47) simplifies to $${1\over 2}\,\sum_{\mu\in
\Gamma_+}\,(\sum_{i\in I^{\hh}}(tr\,ad\,\,b_i\vert
\g^{\mu})b_i)\,a\,u = v^{(1)}\,a\,u\eqno (2.48)$$ But then by (2.44) one has 
$$\eqalign{{1\over 2}\,\sum_{\mu\in
\Gamma_+}\,(\sum_{i\in I^{\hh}}(tr\,ad\,\,b_i\vert
\g^{\mu})b_i)&= \sum_{i\in I^{\hh}}\rho_1(b_i)b_i\cr
&= y_{\rho_1}\cr}$$ Recalling (2.47) we have proved \vskip 1pc {\bf Lemma 2.8.} {\it Let $a\in
C(\p^0)$. Then
$$v^{(1)}\,a\,u = y_{\rho_1}\,a\,u$$} \vskip 1pc Recall that $\square = \square' +
\square''$ (see (2.6) and (2.7)). We now find a condition on $\lambda$ to insure that
$v_{\lambda,a}$ is a Dirac cocyle (assuming that $a\in C(\hh_{\p})$). 
\vskip 1pc {\bf Theorem 2.9.} {\it Let $\lambda \in \hh^*$ and let $a\in C(\hh_{\p})$. 
 Recall $v_{\lambda,a} = v_{\lambda} \otimes a\,u\in V_{\lambda}\otimes L$. Then
$$\xi_{\lambda}(\square'')v_{\lambda,a} = v_{\lambda}\otimes y_{\rho}\,a\,u\eqno (2.49)$$
Furthermore $$\xi_{\lambda}(\square)v_{\lambda,a} = v_{\lambda}\otimes
y_{\lambda +\rho}\,a\,u\eqno (2.50)$$ In particular $$\xi_{\lambda}(\square)v_{\lambda,a} =
0\eqno (2.51)$$ in case $(\lambda +\rho)\vert \hh_{\p}= 0$.} \vskip 1pc {\bf Proof.} Equation
(2.49) follows from (2.27), (2.33), Lemma 2.8, (2.39) and (2.40). But then (2.50) follows
from Lemma 2.3 and (2.48). The equation (2.51) is immediate from the definition of $y_{\lambda
+\rho}$ (see (2.11)). QED

\eject
 \centerline{\bf 3. Non-vanishing Dirac
cohomology}\vskip 1.5pc 3.1. Henceforth we assume that $\lambda\in \hh^*$ is an arbitrary
element satisfying $\break$ $(\lambda +\rho)\vert \hh_{\p}= 0$. We will establish that
$v_{\lambda,a}$, for any
$0\neq a\in C(\hh_{\p})$ defines a nonzero Dirac cohomology class. 

Let $\n^0_-$ be the span of the root vectors $e_{-\varphi}$ for $\varphi\in \Delta^0_+$ so that
one has the triangular decomposition $$\p^0 = \hh_{\p} + \n^0 + \n^0_-\eqno (3.1)$$ Put
$\m_{\p}^- = \n^0_- + \p_-$ so that (see (2.2)) one also has the direct sum $$\p = \hh_{\p} +
\m_{\p} + \m_{\p}^-\eqno (3.2)$$ If $\n_-$ is the span of all root vectors $e_{-\varphi}$
for $\varphi\in \Delta_+$ note that $$\m_{\p}^- \s \n_-\eqno (3.3)$$ since $\p_-$ is spanned
by eigenvectors of $ad\,f_{\r}$ for negative eigenvalues. For any subspace
$\a\s
\p$ let
$C(\a)\s C(\p)$ be the Clifford algebra (with respect to $B_{\p}$) generated by $\a$ (and of
course 1). Clearly 
$C(\a)=\wedge\a$. Since $\m_{\p}^-$ is obviously isotropic it follows that exterior and
Clifford multiplication are the same in $C(\m_{\p}^-)$. Now by (2.3), (2.4) and (3.2) note
that the map $$C(\hh_{\r} + \m_{\p}^-) \to L,\qquad w\mapsto w\,u\eqno (3.4)$$ is a linear
isomorphism. Let $C_*(\m_{\p}^-)$ be the ideal in $C(\m_{\p}^-)$ generated by $\m_{\p}^-$. It
then follows from (3.4) that one has a direct sum $$L = C(\hh_{\p})u \oplus
C_*(\m_{\p}^-)C(\hh_{\p})u \eqno (3.5)$$ Now if $\lambda',\lambda''\in \hh^*$ we will say that
$\lambda''$ is less than $\lambda'$ (or $\lambda'$ is greater than $\lambda''$) and write
$\lambda' > \lambda''$ in case $\lambda'-\lambda''$ is a nontrivial sum of positive roots.
Let $V_{\lambda,*}$ be the span of all weight vectors, of some weight $\lambda'$ in
$V_{\lambda}$,  where $\lambda > \lambda'$. Then clearly one has $$V_{\lambda} = \Bbb
C\,v_{\lambda} \oplus V_{\lambda,*}\eqno (3.6)$$ For notational convenience let $M\s
V_{\lambda}\otimes L$ be defined by putting $$M = (V_{\lambda}\otimes
C_*(\m_{\p}^-)C(\hh_{\p})u)\,\, + \,\, V_{\lambda,*}\otimes L\eqno (3.7)$$ It then follows
from (3.5) and (3.6) that
$$\eqalign{V_{\lambda,L}& = V_{\lambda}\otimes L\cr &= (\Bbb C\,v_{\lambda}\otimes
C(\hh^{\p})u)\oplus M\cr}\eqno (3.8)$$ 

Let $H_{D}(V_{\lambda,L})$ be the
Dirac cohomology defined by $\xi_{\lambda}(\square)$. See (0.3). By Theorem (2.9) the map from
a cocyle to the  corresponding cohomology class defines a linear map $$\Bbb
C\,v_{\lambda}\otimes C(\hh^{\p})u\to H_{D}(V_{\lambda,L})\eqno (3.9)$$ We will show that
(3.9) is injective.
\vskip 1pc {\bf Proposition 3.1}. {\it To show that (3.9) is injective it suffices to prove
that $M$ is stable under the action of $\xi_{\lambda}(\square)$.} \vskip 1pc {\bf
Proof.} This is immediate from (3.8) since $\Bbb C\,v_{\lambda}\otimes
C(\hh^{\p})u \s Ker\,\xi_{\lambda}(\square)$ by Theorem 2.9. QED \vskip 1pc 

	3.2. To show that $M$ is stable under $\xi_{\lambda}(\square)$ we first establish\vskip 1pc
{\bf Lemma 3.2.} {\it The space
$M$ is stable under
$\xi_{\lambda}(\square')$.} \vskip 1pc {\bf Proof.} We use the notation of \S 2.2 where
$\square'$ is given by (2.12) and the basis $b_i,\,i\in I$, is defined as in \S 2.2. To prove
the lemma it suffices to show that $M$ is stable under $\xi_{\lambda}(b_i\otimes d_i)$ for
any
$i\in I$. Assume first that $i\in I^0\cup I'$. It is obvious from (2.18) that
$$\hbox{either}\,\,b_i\in \m_{\p}^-\,\,\hbox{or}\,\,d_i\in \m_{\p}^-\eqno (3.10)$$ If $b_i\in
\m_{\p}^-$ then $b_i\in \n_-$ by (3.3) so that clearly $$\eqalign{\xi_{\lambda}(b_i\otimes
d_i)(M)&\s \xi_{\lambda}(b_i\otimes
d_i)(V_{\lambda}\otimes L)\cr
&\s V_{\lambda,*}\otimes L\cr
&\s M\cr}$$ If $d_i\in \m_{\p}^-$ then by (3.4) $$\eqalign{\xi_{\lambda}(b_i\otimes
d_i)(M)&\s \xi_{\lambda}(b_i\otimes
d_i)(V_{\lambda}\otimes L)\cr &\s V_{\lambda}\otimes  C_*(\m_{\p}^-)C(\hh_{\p})u\cr
&\s M\cr}$$ Now assume $i\in I^{\hh}$ so that $b_i\otimes d_i= b_i\otimes b_i$ where
$b_i\in\hh_{\p}$ by (2.17). But obviously $$\eqalign{b_i\,V_{\lambda,*}&\s
V_{\lambda,*}\cr  b_i\,\,C_*(\m_{\p}^-)C(\hh_{\p})u &\s
C_*(\m_{\p}^-)C(\hh_{\p})u\cr}$$ so that
$M$ is stable under $\xi_{\lambda}(b_i\otimes d_i)$ in this case as well. QED \vskip 1pc
3.3. It remains only to show that $M$ is stable under $\xi_{\lambda}(\square'')$. But now for
$z\in V_{\lambda}$ it is obvious that
$z\otimes L$ is stable under $\xi_{\lambda}(\square'')$. The question is then reduced to
considering only $L$. In fact one immediately has \vskip 1pc {\bf Lemma 3.3.} {\it If
$$v\,C_*(\m_{\p}^-)C(\hh_{\p})u\s C_*(\m_{\p}^-)C(\hh_{\p})u\eqno (3.11)$$ then $M$ is stable
under $\xi_{\lambda}(\square'')$.} \vskip 1pc We now proceed to establish the inclusion
(3.11). 

Let $SO(\p)$ be the special orthogonal group with respect to $B_{\p}$. One has a homomorphism
$$\nu:\r \to Lie\,SO(\p)\eqno (3.12)$$ where if $y\in \p$ and $x\in \r$ then
$\nu(x)y = [x,y]$. Now $\wedge^2\p\s C(p)$ has the structure of a Lie algebra under Clifford
commutation and one has a Lie algebra isomorphism $\wedge^2\p\cong Lie\,SO(\p)$. See (1.7) in
[K1]. Furthermore there exists a Lie algebra homomorphism $\nu_*:\r \to \wedge^2\p$ so
that for $x\in \r$ and $y\in \p$ one has, using commutation in $C(\p)$,
$$\eqalign{\nu(x)y &= [\nu_*(x),y]\cr
&=-2\,\iota(y)\nu_*(x)\cr}\eqno (3.13)$$ where $\iota(y)$ is the operator of interior
product   of $\wedge\,\p$ by $y$. See (1.8) and (1.11) in [K1]. Let $\{b_i\}$ and $\{d_i\}$,
for $i\in
I$, be the basis and dual basis of $\p$ defined as in \S 2.1. Let $I^+ = \{i\in I\mid b_i\in
\p^+\}$. See (1.14). Note $\{b_i\}, \,i\in I^+$, is a basis of $\p^+$ and $\{d_i\}, \,i\in
I^+$, is a basis of $\p^-$ (see (2.15)) so that taken together $\{b_i\}\cup\{d_i\},\,i\in
I+,$ is a basis of
$\p'$. \vskip 1pc {\bf Proposition 3.4.} {\it For any
$x\in
\hh_{\r}$ one has, using the notation of (2.26),
$$\nu_*(x) = {1\over 2}\,\sum_{i\in I^+} \gamma_i(x)\,b_i\wedge d_i\eqno (3.14)$$}\vskip 1pc
{\bf Proof.} Let $w\in \wedge^2\p$ be given by the right side of (3.8). Then clearly, for
$i\in I^+$, $-2\iota(d_i)w= -\gamma_i(x)d_i$ and $-2\iota(b_i)w = \gamma_i(x)b_i$. On the
other hand
$-2\iota(y)w = 0 $ if $y\in\p^0$. But, by (3.13), these same equations are satisfied if
$\nu_*(x)$ replaces
$w$. This proves $w= \nu_*(x)$. QED \vskip 1pc Let $\Lambda\s \hh_{\r}^*$ be
the real space of all (complex) linear functionals $\beta$ on $\hh^{\r}$ such that there exists
$\gamma\in \hh^*$ with the property that (1) $\gamma\vert\hh_{\r} = \beta$ and (2)
$\gamma(\hh^{\#})\s \Bbb R$. It is immediate that any $\beta\in \Lambda$ extends uniquely,
as a linear functional, on the complex subspace of $\hh$ spanned by $\kappa(\hh_{\r})$. (See
\S 1.1). In particular $\beta(f_{\r})$ is well-defined for $\beta\in \Lambda$. Obviously
(2.45) implies that 
$\Gamma\s \Lambda$. It follows therefore that $\rho_{\p}\in \Lambda$ where for any $x\in
\hh_{\r}$ $$\eqalign{\rho_{\p}(x) &= {1\over 2}\, tr\,\, \nu(x)\vert \p^+\cr &= {1\over
2}\,\sum_{i\in I^+}\gamma_i(x)\cr}\eqno (3.15)$$\vskip .5pc {\bf Proposition 3.5.} {\it Let
$x\in \hh_{\r}$. Then $$\nu_*(x)\,u = \rho_{\p}(x)\,u\eqno (3.16)$$}\vskip 1pc {\bf Proof.}
Let $i\in I^+$. Then $b_i\wedge d_i = - d_i\,b_i +1 $ by (2.32). But $b_i\,u = 0$ by (2.3).
But then (3.16) follows from (3.14) and (3.15). QED \vskip 1pc Now consider the action of
$\hh_{\r}$ on $L$ defined, for $x\in \hh_{\r}$, by left multiplication on $L$ by $\nu_*(x)$.
\vskip 1pc {\bf Lemma 3.6}. {\it If $s\in L$ is an $\hh_{\r}$-weight vector with weight
$\beta\in \hh_{\r}^*$ and $i\in I$, then $b_i\,s_i$ is an $\hh_{\r}$-weight vector with weight
$\gamma_i +\beta$.}\vskip 1pc {\bf Proof.} Let $x\in \hh_{\r}$. Then $[\nu_*(x),b_i]=
\gamma_i(x)\,b_i$ by (3.13). But $$\eqalign{\nu_*(x)\,b_i\,s &= [\nu_*(x),b_i]\,s +
b_i\,\nu_*(x)s\cr&= (\gamma_i(x) + \beta(x))b_i\,s \cr}$$ QED \vskip 1pc It is immediate
from (3.4), Proposition 3.5 and Lemma 3.6 that if $\beta$ is an $\hh_{\r}$-weight in $L$ then
$\beta\in \Lambda$. Let $\Lambda_L\s \Lambda$ be the set of all $\hh_{\r}$-weights in $L$ and
for any $\beta\in \Lambda_L$ let $L^{\beta}$ be the corresponding weight space so that one has
the direct sum $$L = \sum_{\beta\in \Lambda_L} L^{\beta}\eqno (3.17)$$ \vskip .5pc
{\bf Proposition 3.7.} {\it For any $\beta\in \Lambda_L$ the weight space $L^{\beta}$ is stable
under left multiplication by the element $v\in \wedge^3\p$.} \vskip 1pc {\bf Proof.} It is
immediate from (2.8) that $v$ is invariant under the representation $\theta_{\nu}$ of
$\hh_{\r}$ on $\wedge\p$ using the notation of (1.12) in [K1]. But by (1.12) in [K1] this
implies that, with respect to Clifford multiplication,
$v$ commutes with $\nu_*(x)$ for all $x\in \hh_{\r}$. The proposition then follows
immediately. QED \vskip 1pc For notational convenience let $L_- = C_*(\m_{\p}^-)C(\hh_{\p})u$
so that by (3.4) one has the direct sum $$ L= C(\hh_{\p})\,u \oplus L_-\eqno (3.18)$$ Our
problem is to show that $L_-$ is stable under left multiplication by $v$. Let $C_*(\n^0_-)$ be
the ideal in $C(\n^0_-)$ generated by $\n^0_-$ and let $C_*(\p_-)$ be the ideal in $C(\p_-)$
generated by $p_-$. Note that $L_-$ can be written $$L_- = (C_*(\n^0_-)C(\hh_{\p})u) \oplus
(C_*(\p_-) C(\n^0_-)C(\hh_{\p})u)\eqno (3.19)$$ 
Now for $\beta,\beta'\in \Lambda$ we will write
$\beta >
\beta'$ (and say that $\beta$ is higher than $\beta'$) in case
$\beta(f_{\r})>\beta'(f_{\r})$. It is then immediate from Proposition 3.5, Lemma 3.6, (3.18)
and (3.19) that $\rho_{\p}$ is the highest $\hh_{\r}$ weight in $L$
and $$\eqalign{C(\n_0)C(\hh_{\r})u &= L^{\rho_{\p}}\cr
C_*(\p_-) C(\n^0_-)C(\hh_{\p})u&= \sum_{\beta\in
\Lambda_L,\,\,\rho_{\p}>\beta}L^{\beta}\cr}\eqno (3.20)$$\vskip .5pc We can now simplify our
problem. \vskip 1pc {\bf Proposition 3.8.} {\it To prove that $v\,L_-\s L_-$ it suffices only
to show that the subspace $C_*(\n^0_-)C(\hh_{\p})u$ is stable under left multiplication by
$v$.}
\vskip 1pc {\bf Proof.} By (3.20) and Proposition 3.7 it follows that $C_*(\p_-)
C(\n^0_-)C(\hh_{\p})u$ is stable under left multiplication by $v$. But then Proposition 3.8
follows from (3.19). QED\vskip 1pc Now recall (see (2.27)) we have written $v = v^{(0)} +
v^{(1)}$. \vskip 1pc {\bf Proposition 3.9.} {\it The space $C_*(\n^0_-)C(\hh_{\p})u$ is
stable under left multiplication by $v^{(1)}$. } \vskip 1pc {\bf Proof.} Since
$y_{\rho_1}\in \hh_{\p}$ the proof follows immediately from Lemma 2.8. QED \vskip 1pc We are
reduced finally to showing that
$C_*(\n^0_-)C(\hh_{\p})u$ is stable under left multiplication by $v^{(0)}$. Let $L_0 =
C(\p^0)u$ so that $L_0$ is a cyclic $C(\p^0)$-module under left multiplication. On the other
hand recalling the triangular decomposition (3.1) and recalling the definition of $u = u_0u_+$
in \S 2.1 one has $$L_0 = C(\hh_{\p})u \oplus C_*(\n^0_-)C(\hh_{\p})u\eqno (3.21)$$ \vskip
.5pc {\bf Remark 3.10.} Note that $L_0$ is stable under left multiplication by $v^{(0)}$
since, by definition, $v^{(0)}\in C(\p^0)$. \vskip 1pc 3.4. Recalling that $\p^0$ is a
reductive Lie algebra and $B_{\p}\vert \p^0$ is nonsingular let $\sigma:\hh_{\p}\to Lie\,
SO(\p^0)$ be defined so that for $y\in \hh_{\p}$ and $z\in \p^0$ one has $\sigma(y)z =
[y,z]$. Going back again to \S 1.5 in [K1] one has a Lie algebra homomorphism
$\sigma_*:\hh_{\p}\to
\wedge^2\p^0$ so that for
$y\in \hh_{\p}$ and $z\in \p^0$ one has $\sigma(y)z = [\sigma_*(y),z]$. Noting that
$\sigma(y)z = 0$ for $z\in\hh_{\p}$ the argument establishing (3.14) readily establishes
\vskip 1pc {\bf Proposition 3.11.} {\it For any $y\in \hh_{\p}$ one has $$\sigma_*(y) =
{1\over 2}\,\sum_{\varphi\in \Delta^0_+}\varphi(y)\,e_{\varphi}\wedge e_{-\varphi}\eqno
(3.22)$$}\vskip 1pc Now, recalling the definition of $\rho_0$ in \S 2.2 (on the line
following (2.28)), one has 
$\rho_0  =  {1\over 2} \sum_{\varphi\in \Delta^0_+}\varphi$. The argument establishing
Proposition 3.5 yields \vskip 1pc {\bf Proposition 3.12.} {\it Let $y\in \hh_{\p}$. Then
$$\sigma_*(y)\,u = \rho_0(y)\,u\eqno (3.23)$$} \vskip 1pc The nonsingularity of
$B_{\p}\vert\p^0$ implies the nonsingularity of $B_{\p}\vert [\p^0,\p^0]$. But then,
recalling (1.11), one has the nonsingularity of $B_{\p}\vert \d_{\p}$ since $\d_{\p}$ is a
Cartan subalgebra of the semisimple Lie algebra $[\p^0,\p^0]$. Recalling (1.8) one has
$\d_{\p}\s\hh_{\p}$. Let
$\ee_{\p}$ be the $B_{\g}$-orthocomplement of $\d_{\p}$ in $\hh$ so that $$\hh = \d_{\p} +
\ee_{\p}\eqno (3.24)$$ is a $B_{\g}$-orthogonal direct sum. Let $\varphi\in \Delta^+$. Then
one must have
$[e_{\varphi},e_{-\varphi}]\in \d_{\p}$ so that $$\varphi\vert\ee_{\p} = 0\eqno (3.25)$$
With respect to the decomposition (3.19) let
$f_{\p}\in \d_{\p}$ be the component in $\d_{\p}$ of the regular hyperbolic element $f\in
\hh^{\#}$. But then for any $\varphi\in \Delta^0$ one has $\varphi(f_{\p}) >0$ or
$\varphi(f_{\p})<0$ according as $\varphi\in \Delta^0_+$ or $\varphi\in \Delta^0_-$. If
$\delta,\delta'\in
\hh_{\p}^*$ we will say that $\delta$ is higher than $\delta'$ and write $\delta >\delta'$ if
$(\delta-\delta'')(f_{\p})$ is a positive real number. Now let $D\s \hh_{\p}^*$ be the set of
weights for the action of $\hh_{\p}$ on $L_0$ where $y\in \hh_{\p}$ operates as left
multiplication by $\sigma_*(y)$. For any $\delta\in D$ let $L_0^{\delta}$ be the weight space
for the weight $\delta$. If
$\varphi\in
\Delta^0$ and
$t\in L_0$ is a weight vector with weight $\delta$, then the argument establishing Lemma 3.6
also establishes that $e_{\varphi}\,t$ is a weight vector with weight $\tilde{\varphi} +
\delta$ where
$\tilde{\varphi}= \varphi\vert \hh_{\p}$. But then Proposition 3.12 and (3.21) imply \vskip
1pc {\bf Proposition 3.13.} {\it Let
$\tilde{\rho_0} = \rho_0\vert \hh_{\p}$. Then $\tilde{\rho_0}\in D$ and $\tilde{\rho_0}$ is the
highest weight. Moreover
$$\eqalign{C(\hh_{\p})u &= L_0^{\tilde{\rho_0}}\cr
C_*(n_-^0)C(\hh_{\p})u &= \sum_{\delta\in D,\,\tilde{\rho_0}>\delta} L_0^{\delta}\cr}\eqno
(3.26)$$}\vskip 1pc 3.5. We can establish the final step. \vskip 1pc {\bf Proposition 3.14.}
{\it The space $C_*(n_-^0)C(\hh_{\p})u $ is stable under left multiplication by $v^{(0)}$.}
\vskip 1pc {\bf Proof.} It is clear from the definition of $v^{(1)}$ (see (2.27) that
$$(v^{(1)},y\wedge y'\wedge y'') = 0$$ for any $y,y',y\in \p^0$. But $v^{(0)}\in \wedge^3\p^0$
and hence $$([y,y'],y'') = -2 (v^{(0)},y\wedge y'\wedge y'')\eqno(3.27)$$ for any $y,y',y\in
\p^0$ by (2.8). But then it follows immediately from (3.27) that $v^{(0)}$ is invariant under
$\theta_{\sigma}(z)$ for any $z\in \hh_{\p}$ using the notation of (1.12) in [K1]. But then
$v^{(0)}$ commutes with $\sigma_*(z)$ in $C(\p^0)$ for any $z\in \hh_{\p}$ by (1.12) in
[K1].
 It follows therefore that any weight space $L_0^{\delta}$ is stable under left
multiplication by $v^{(0)}$. But by Proposition 3.13 this implies that
$C_*(n_-^0)C(\hh_{\p})u $ is stable under left multiplication by $v^{(0)}$. QED \vskip 1pc
We have proved \vskip 1pc {\bf Theorem 3.15.} {\it If $\lambda\in \hh^*$ is such that
$\lambda +\rho$ vanishes on $\hh_{\p}$ then
$V_{\lambda,L} = V_{\lambda}\otimes L$ has nonvanishing Dirac cohomology. In fact the map
(3.9) $$\Bbb C v_{\lambda}\otimes C(\hh_{\p})u\to H_{D}(V_{\lambda,L})$$ is injective. One
notes that $dim\,C(\hh_{\p})u = 2^k$ where $k= dim\,\hh_{\p}$.} \vskip 1.5pc \centerline{\bf
4. Consequences of Theorem 3.15}\vskip 1.5pc 4.1.  The homomorphism $\nu_*:\r \to
\wedge^2\p\s C(\p)$ (see (3.12) and (3.13)) defines a homomorphism $$\zeta:U(\r)\to
U(\g)\otimes C(\p)\eqno (4.1)$$ where if $x\in \r$ then $\zeta(x) = x\otimes 1 + 1\otimes
\nu_*(x)$. See \S 2.15 in [K1]. This defines the structure of an $\r$-module on $U(\g)\otimes
C(\p)$. Let $(U(\g)\otimes C(\p))^{\r}$ denote the algebra of $\r$-invariants in $U(\g)\otimes
C(\p)$. Let $Z(\g )$ and
$Z(\r)$, respectively, be the centers of
$U(\g)$ and $U(\r)$. One notes that $\square\in (U(\g)\otimes C(\p))^{\r}$ and also
$Z(\g)\otimes 1$ and $\zeta(Z(\r))$ are subalgebras of $(U(\g)\otimes C(\p))^{\r}$.
In case $\r$ is symmetric the cubic term in $\square$ vanishes. The main result in [HP]
(Theorem 3.4) is a statement for the case where
$\r$ is symmetric. However, as noted in the
Appendix, the proof in [HP] is valid in the general case considered here (i.e., the case
where $\r$ is arbitrary, subject only to the condition that
$\r$ is reductive and
$B_{\g}\vert
\r$ is nonsingular and the cubic Dirac operator $\square$ replaces the more familiar Dirac
operator in [HP] and [P]). In addition Corollary 3.5 in [HP] is also valid in the
general case considered here. That is, one has a unique map 
$$\eta_{\r}:Z(\g)\to Z(\r)\eqno (4.2)$$ with the property that, for $p\in Z(\g)$, there
exists
$\omega\in  (U(\g)\otimes C^{odd}(\p))^{\r}$ such that $$p\otimes 1-
\zeta(\eta_{\r}(p))=
\square \,\omega + \omega\, \square\eqno (4.3)$$ Furthermore $\eta_{\r}$ is an algebra
homomorphism. But
$H_D(V_{\lambda,L})$ is a module for both
$Z(\g)\otimes 1$ and $\zeta(Z(\r))$. But (4.3) implies that for any $p\in Z(\g)$ one has
$$p\otimes 1 = \zeta(\eta_{\r}(p))\,\,\hbox{on}\,\,H_D(V_{\lambda,L})\eqno (4.4)$$ Let
$W_{\g}$ be the Weyl group of the pair $(\hh,\g)$ operating on the symmetric algebra $S(\hh)$
over $\hh$ and let $W_{\r}$ be the Weyl group of the pair $(\hh_{\r},\r)$ operating on the
symmetric algebra $S(\hh_{\r})$ over $\r$. One has the Harish-Chandra algebra isomorphisms
$$\eqalign{A_{\g}&:Z(\g)\to S(\hh)^{W_{\g}}\cr
A_{\r}&:Z(\r)\to S(\hh_{\r})^{W_{\r}}\cr}\eqno (4.5)$$ Let $\rho_{\r}\in \hh_{\r}^*$ be
defined so that for $x\in \hh_{\r}$ one has $\rho_{\r}(x) = {1\over 2}\,tr\,\,
ad\,x\vert\n_{\r}$ (see 1.18)). It follows immediately from (1.17) that on
$\hh_{\r}^*$, $$\rho_{\p} + \rho_{\r} = \rho\vert\hh_{\r}\eqno (4.6)$$ Let $E_{\lambda,L}\s
H_D(V_{\lambda,L})$ be the image of (3.9). On the other hand the elements of $\Bbb
C\,v_{\lambda}\otimes C(\hh_{\r})u$ are $\g\otimes 1$ highest weight vectors with highest
weight
$\lambda$ and these elements are highest weight vectors for $\zeta(\r)$ with highest weight
$\lambda\vert \hh_{\r} + \rho_{\p}$. But this establishes the following generalization of
Theorem 5.5 in [HP] (case where $\r$ is symmetric). 
\vskip 1pc {\bf Theorem 4.1.} {\it Let
$p\in Z(\g)$ and $q\in Z(\r)$. Then $p\otimes 1$ reduces to the scalar
 $A_{\g}(p)(\lambda+ \rho)$ on $E_{\lambda,L}$ and (by 4.6) $\zeta(q)$ reduces to the scalar
$A_{\r}(q)((\lambda + \rho)\vert\hh_{\r})$ on $E_{\lambda,L}$.} \vskip 1pc But now since
$(\lambda+\rho)\vert \hh_{\p} = 0$ and there is no restriction on $\lambda\vert \hh_{\r}$,
this completely determines the map $\eta_{\r}$. In fact let $\phi_o:\hh\to \hh_{\r}$ be the
projection relative to the decomposition $\hh = \hh_{\r} + \hh_{\p}$. Then $\phi_o$
extends to a homomorphism $S(\hh)\to S(\hh_{\r})$ and clearly induces a homomorphism
$$\phi:S(\hh)^{W_{\g}}\to S(\hh_{\r})^{W_{\r}}\eqno (4.7)$$ Given the Harish-Chandra
isomorphisms $A_{\g}$ and $A_{\r}$ the map $\eta_{\r}$ is given by completing a commutative
diagram. Since one must have $A_{\g}(p)(\lambda+ \rho)=
A_{\r}(\eta_{\r}(p))((\lambda +
\rho)\vert \hh_{\r})$ for all $p\in Z(\g)$ and all $\lambda\in\hh^*$ such that
$(\lambda+\rho)\vert \hh_{\p}=0$ we have established the following generalization of
Theorem 5.5 in [HP]. It follows from the observations in \S 5 below that it is also a
generalization of Proposition 3.43, (5.18) and (5.19) in [K1]. 
\vskip 1pc {\bf Theorem 4.2.} {\it The map
$\eta_{\r}:Z(\g)\to Z(\r)$ is uniquely determined so that the following diagram is
commutative
$$
\def\mapright#1{\smash{
\mathop{\longrightarrow}\limits^{#1}}}
\def\mapdown#1{\Big\downarrow
\rlap{$\vcenter{\hbox{$\scriptstyle#1$}}$}}
\matrix
{Z(\g)&\mapright {\eta_{\r}}\, &Z(\r)\cr
\mapdown {A_{\g}}&&\mapdown {A_{\r}}\cr
S(\hh)^{W_{\g}}&\mapright{\phi}&S(\hh_{\r})^{W_{\r}}\cr}$$}

The map $\phi$ is well known in the theory of the cohomology of compact homogeneous
spaces. Actually what is utilized in that theory is the map $S(\hh^*)^{W_{\g}}\to
S(\hh_{\r}^*)^{W_{\r}}$ induced by restriction of functions. However this is same as
$\phi$ if $\hh$ and $\hh^*$ are identified and $\hh_{\r}$ and $\hh_{\r}^*$ are
identified using $B_{\g}$. Assume $G$ is a compact connected semisimple Lie group and
$\g$ is the complexification of $Lie\,G$. Let $R\s G$ be any connected compact subgroup and
let
$\r$ be the complexification of $Lie\,R$. Obviously we can choose $B_{\g}$ so that
$B_{\g}\vert \r$ is nonsingular (e.g., let $B_{\g}$ be the Killing form). The 
map $\eta_{\r}$ induces the structure of a $Z(\g)$-module on $Z(\r)$. On the other
hand the infinitesimal character for the module $V_{\lambda}$ when $\lambda =-\rho$
defines the structure of a $Z(\g)$-module on $\Bbb C$. As a consequence of a
well-known theorem of H. Cartan (see \S 9 in [C]) one has \vskip 1pc {\bf Theorem 4.3.} {\it
There exists an isomorphism $$H^*(G/R,\Bbb C) \cong Tor_*^{Z(\g)}(\Bbb C, Z(\r))\eqno
(4.8)$$}\vskip 1.5pc \centerline{\bf 5. The case where $rank\,\r = rank\,\g$ and
$dim\,V_{\lambda}<\infty$}\vskip 1.5pc 5.1. Let the notation be as \S0.1 so that $V$ is an
arbitrary $\g$-module. It is clear that $\zeta(\r)$ commutes with $\square$ so that
the Dirac cohomology, $H_D(V\otimes L)$ has the structure of an $\r$-module. Of course
$Ker\,\square_{V}\s Ker\,\square_{V}^2$. Note that the special case $$Ker\,\square_{V}=
Ker\,\square_{V}^2\eqno (5.1)$$ occurs if and only if $$Ker\,\square_{V}\cap
Im\,\square_{V} = 0\eqno (5.2)$$ If (5.1), or equivalently (5.2), occurs then we may regard
$H_D(V\otimes L)\s V\otimes L$ where in fact one has $$\eqalign{H_D(V\otimes
L)&=Ker\,\square_{V}\cr &=Ker\,\square_{V}^2\cr}\eqno (5.3)$$ In this section we would like
to formulate results in [K1] and [K2], especially results beginning with
\S 3 in [K1], in terms of Dirac cohomology. Assume, as in \S 3 of [K1], that $rank\,r
=rank\,\g$ so that $\hh =\hh_{\r}$ and $\hh_{\p} = 0$. Note that, in this case, the
restriction on $\lambda$ in Theorem 3.15 disappears. Also in this case 
$$S= L\eqno (5.4)$$ where $S$ is the $C(\p)$-spin module of \S 3.1 in [K1]. See (3.11) in
[K1]. Next assume that
$\lambda$ is dominant with respect to $\b$ and integral with respect to $\g$. But then
$V_{\lambda}$ is finite dimensional and $\g$-irreducible. Consider $H_D(V_{\lambda}\otimes
S)$. Using the notation of (4.8) let $R\s G$ be any connected compact subgroup having the
same rank as
$G$. Up to conjugacy we can take $R$ to be defined as in
\S 5.21 in [K1] so that $\r$ is the complexification of $Lie\,R$. In this section, as in
[K1], let $d$ be the Euler characteristic of $G/R$. We have written $W$ for the Weyl group
$W_{\g}$ in [K1]. One has $W_{\r}\s W$ and one knows that the index of $W_{\r}$ in $W$ is
$d$. See e.g., (5.32) in [K1]. Let $W^1\s W$ be the set of representatives of the right
cosets of
$W_{\r}$ in $W$ defined as in (3.24) in [K1] so that $d = card\,\,W^1$. For any $\tau\in
W^1$ let $\tau\bullet\lambda = \tau(\lambda + \rho)-\rho_{\r}$. Then $\tau\bullet\lambda$
is dominant with respect to $\b_{\r}$ and integral for the simply-connected covering group
of $R$. In particular if $Z_{\tau\bullet\lambda}$ is an irreducible $\r$-module with
highest weight $\tau\bullet\lambda$ then $Z_{\tau\bullet\lambda}$ is finite dimensional.
Also $Z_{\tau_i\bullet\lambda},\,i=1,2,$ are inequivalent for $\tau_i\in W^1$ where
$\tau_1\neq \tau_2$. See \S 3.22 in [K1]. But now Theorems 4.17 and 4.24
in [K1] imply
\vskip 1pc {\bf Theorem 5.1.} {\it Assume $rank\,\r=rank\,\g$ and $\lambda\in \hh^*$ is
dominant and integral with respect to $G$. Then $Z_{\tau\bullet\lambda}$ occurs with
multiplicity one in $V_{\lambda}\otimes S$, for any $\tau\in W^1$, so that we can
unambiguously regard $Z_{\tau\bullet\lambda}\s V_{\lambda}\otimes S$. Furthermore the
condition (5.1) is satisfied and (recalling (5.3))
$$H_D(V_{\lambda}\otimes S) = \sum_{\tau\in W^1} Z_{\tau\bullet\lambda}\eqno (5.5)$$ In
particular $H_D(V_{\lambda}\otimes S)$, as an $\r$-module, is multiplicity-free and
decomposes into a sum of
$d$ irreducible components, where $d$ is the Euler number of $G/R$.} \vskip 1pc {\bf Remark
5.2.} In the case where $\r$ is the Levi factor of a parabolic subalgebra of $\g$ we have
shown in [K2] that Theorems 4.17 and 4.24 in [K1] imply the Bott-Borel-Weil theorem (BBW).
This may formulated in terms of Dirac cohomology. In case $\r$ is the Levi factor of a
parabolic subalgebra of $\g$ the argument in [K2] shows that BBW is a consequence
of Theorem 5.1 together with the construction of
$Z_{\tau\bullet\lambda}$ given in Theorem 4.17 of [K1]. \vskip 1pc 5.2. As mentioned above,
Theorem 4.2, for the case where $rank\,\r = rank\,\g$ appears in [K1]. See \S 5, especially
equations (5.18) and (5.19), in [K1]. In more detail, the map $\phi$ in the present
case, is injection so that the map
$\eta_{\r}$ is injective. The image of $Z(\g)$ in $Z(\r)$ has been denoted by
$Z_{\g}(\r)$ in [K1]. Let the notation be as in Theorem 5.1. The set
$\{Z_{\tau\bullet\lambda}\},\,\,\tau\in W^1$, of representations of $\r$ is referred to in
[K1] as a multiplet. Recalling (5.5), a verification of equation (4.4) 
is the statement that the infinitesmal character of $Z(\r)$ for all the members of a
multiplet remains the same when restricted to $Z_{\g}(\r)$ and that, furthermore, the
restriction is given by the infinitesmal character of $Z(\g)$ for the $\g$-representation 
$V_{\lambda}$. But this and more is stated in Proposition 3.43 of [K1] together with (5.18)
and (5.19). \vskip 1pc {\bf Remark 5.3.} In a certain sense matters have come full circle.
Consider the case where $\g$ is of type
$F_{4}$ and $\r$ is of type $B_4$ (i.e., $R\cong Spin\,9$). In that case $d=3$ so the
multiplets are triplets. That the members of each triplet had remarkable properties in
common was the empirical discovery of the physicists Ramond and Pengpan. This discovery
inspired the paper [GKRS] which in turn led to [K1]. One of the properties discovered by
Ramond and Pengpan, in the terminology above, is, in retrospect, the statement that
$Z_{\g}(\r)$ operates the same way on each member of any triplet. (We use the term ``in
retrospect" since Ramond and Pengpan were dealing only with $\r$ and were unaware of the
role of
$\g$.) But, with the notion of Dirac cohomology, this behavior of $Z_{\g}(\r)$ is necessarily
the case since (see (5.5))
$H_D(V_{\lambda}\otimes S)$ is just the sum of the members of that triplet which corresponds
to $\lambda$.\vskip 1.5pc
	
                          \centerline{\bf Appendix}\vskip 1.5pc
A.1. One of the properties of $\square$ used in [HP] to establish the main
theorem, Theorem 3.4 in [HP], was Lemma 3.1 in [HP]. Let $Cas_{\g}\in Z(\g)$ be
the $\g$-Casimir element with respect to $B_{\g}$ and let $Cas_{\r}\in Z(\r)$
be the $\r$-Casimir element with respect to $B_{\r}$. Recall that $\r$ is
assumed to be symmetric in [HP]. Lemma 3.1 in [HP] asserts 
$$\square^2 = Cas_{\g}\otimes 1 - \zeta(Cas_{\r}) + constant \eqno (A.1)$$ It
should be noted that the definition of Dirac operator in [HP] differs from
its definition here and in [K1], in the symmetric case, by a factor of $i$.
This factor clearly plays no significant role in our concerns here. 

The equation (A.1) is used in [HP] to define a $\Bbb Z_2$-graded differential 
complex in $(U(\g)\otimes C(\p))^{\r}$ with $ad\,\,\square$ as the coboundary
operator. Here commutation with $\square$ is taken in the $\Bbb Z_2$-graded
sense. However in the
general case we are considering, where $\square$ is the cubic Dirac operator, we
have established (A.1). See Theorem 2.16 in [K1]. The validity of (A.1)
enables one to define this complex in the general case using the cubic Dirac
operator. Also, as in [HP], $Z(\g)\otimes 1$ and $\zeta(Z(\r))$ are, in the
general case, spaces of cocyles. Theorem 3.4 in [HP] asserts $\zeta(Z(\r))$ is
isomorphic to the cohomology of this complex. But again the same argument yields
the same result in the general case. The idea in [HP] is to replace $U(\g)$ by
the symmetric algebra $S(\g)$ and to replace $ad\,\,\square$ by its symbol. A
computation of the symbol leads to the Koszul complex. The proof then follows from
the acyclicity of the Koszul complex. The reason why this argument works in
the general case is that one obtains the same symbol. This is because the cubic
term has no affect on the symbol. It should be noted that our result, for the case
where
$\r=0$, appears in [AM]. 

The validity of Theorem 3.4 in [HP], for the general case, leads to the map
(4.2), which one easily shows, is a homomorphism of algebras. Theorem 4.1 here
determines the map (4.2) in the general case.\vskip 1.5pc 
\centerline{\bf References}\vskip 1.5pc

\item {[AM]} A. Alekseev and E. Meinrenken, {\it The non-commutative Weil algebra},
Inventiones math. {\bf 139}(2000), 135-172
\item {[C]} H. Cartan, {\it La transgression dans un groupe de Lie et dans un espace fibr\'e
principal}, Colloque de Topologie, C.B.R.M. Bruxelles 57-71(1950)
\item {[GKRS]} B. Gross, B. Kostant, P. Ramond, S. Sternberg, {\it The Weyl character
formula, the half-spin representations, and equal rank subgroups},  PNAS
{\bf 95}(1998), 8441-8442
\item {[HP]} J-S. Huang and P. Pand$\check {\hbox{z}}$i\'c, {\it Dirac cohomology, unitary
representations and a proof of a conjecture of Vogan}, JAMS {\bf 15} (2002), 185-202;
electronic publication with Pand$\check {\hbox{z}}$i\'c, Sept. 6, 2001.
\item {[K1]} B. Kostant, {\it A cubic Dirac operator and the emergence of Euler number
multiplets of representations for equal rank subgroups}, Duke Math. Jour. {\bf 100} (1999),
447-501.
\item {[K2]} B. Kostant, {\it A generalization of the Bott-Borel-Weil theorem and Euler
number multiplets of representations}, Letters in Mathematical Physics, {\bf 52}(2000),
61-78
\item {[P]} R. Parthasarathy, {\it Dirac operator and the Discrete series},
 Ann. of Math., {\bf 96}\break (1972), 1-30
\smallskip
\parindent=30pt
\baselineskip=14pt
\vskip 1.9pc
\vbox to 60pt{\hbox{Bertram Kostant}
      \hbox{Dept. of Math.}
      \hbox{MIT}
      \hbox{Cambridge, MA 02139}}\vskip 1pc

      \noindent E-mail kostant@math.mit.edu

\end